\documentclass{article}

\usepackage{amsmath}
\usepackage{amsfonts}
\usepackage{amsbsy}
\usepackage{amssymb}

\newcommand {\ZZ} {\mathbb Z}
\newcommand {\NN} {\mathbb N}
\newcommand {\RR} {\mathbb R}

\newcommand {\BB} {{\mathcal B}}
\newcommand {\CC} {{\mathcal C}}

\newcommand {\LS}[1] { {\mathcal L}_{f_d}\left({#1}\right) }
\newcommand {\PROBA}[1] {\mathbb P\left( #1 \right)}
\newcommand {\EXPECT}[1] {\mathbb E\left( #1 \right)}

\newcommand {\one} {1 \!\! 1}

\newcommand {\DD} {{\mathcal D}}
\newcommand {\Dom}[1] {\DD\left( #1 \right)}
\newcommand {\LEB}[2] {{\mathbb L}_{#1}\left( #2 \right)}
\newcommand {\Cl}[1] {{\mathcal Cl}_{#1}}

\newcommand {\SSAT}[2] {{\mathcal S}_{#1}^{#2}}

\def\endproof{\hfill $\Box$\newline\newline}  
\def\proof{\par\noindent{\it Proof}. \ignorespaces}

\newtheorem{theorem}{Theorem}[section]

\newtheorem{e-proposition}[theorem]{Proposition}

\newtheorem{e-definition}[theorem]{Definition\rm}
\newtheorem{remark}{\it Remark\/}

\begin{document}
\title{Estimating the probability law of the codelength as a function of the
  approximation error in image compression} 
\author{Fran\c{c}ois Malgouyres\footnotemark[2]}

\maketitle
\renewcommand{\thefootnote}{\fnsymbol{footnote}}

\footnotetext[2]{LAGA/L2TI, Universit\'e Paris 13, 99 avenue
    Jean-Batiste Cl\'ement, 93430 Villetaneuse, France. \\
malgouy@math.univ-paris13.fr\\
http://www.math.univ-paris13.fr/$\sim$malgouy/
}

\renewcommand{\thefootnote}{\arabic{footnote}}

\begin{abstract}
After a recollection on compression through a projection onto a polyhedral
set (which generalizes the compression by coordinates quantization), we
express, in this framework, the probability that an image is coded with $K$
coefficients as an explicit function of the approximation error.

\end{abstract}

\section{Introduction}\label{intro-sec}

In the past twenty years, many image processing tasks have been addressed with
two distinct mathematical tools : Image decomposition in a basis and
optimization. 

The first mathematical approach proved very useful and is supported by
solid theoretical foundations which guarantee its efficiency as long as the
basis is adapted to the information contained in images. Modeling the image
content by appropriate function spaces (of infinite dimension), mathematical
theorems tell us how the coordinates of an image in a given basis behave. As
instances, it is possible to characterize Besov spaces (see \cite{Meyerbook1})
and  the space of bounded variation (it is ``almost characterized'' in
\cite{CohenDahmenDaubechiesDevore}) with wavelet coefficients. As a
consequence of these characterizations, one can obtain performance estimate for
practical algorithms (see Th 9.6, pp. 386, in \cite{mallatbook} and
\cite{cdpx,CohenDaubechiesGuleryusOrchard} for more complex analyses). Image
compression and restoration are the typical applications where such analyses
are meaningful.

The optimization methods which have been considered to solve those practical
problems also proved very efficient (see \cite{rof}, for a very famous
example). However, the theory is not able to assess how well they perform,
given an image model. 

Interestingly, most of the community who was primarily involved in the image
decomposition approaches is now focusing on optimization models (see, for
instance, the work on Basis Pursuit \cite{ChenDonoho} or compressed sensing
\cite{DonohoCompSens}). The main reason for that is probably that optimization
provides a more general framework
(\cite{cdll,MalgouyresIeee02,MalgouyresCompression}).

The framework which seems to allow both a good flexibility for practical
applications (see \cite{ChenDonoho} and other papers on Basis Pursuit) and
good properties for theoretical analysis is the projection onto polyhedron or
polytopes. For theoretical studies, it shares simple geometrical properties
with the usual image decomposition models (see \cite{geomoptim}). This might
allow the derivation of approximation results.

The aim of this paper is to state a rigorous\footnote{The theorem concerning
  compression in \cite{geomoptim} is false. The situation turns out to be
  more complex than we thought at the time it was written.} theorem which
relates, asymptotically as the precision grows, the approximation error and
the number of coefficients which are coded (which we abusively call
codelength, for simplicity). More precisely, when the initial datum is assumed
random in a convex set, we give the probability that the datum is coded by $K$
coefficients, as a function of the approximation error (see theorem
\ref{LeTHM} for details).

This result is given in a framework which generalizes the usual coding of the
quantized coefficients (``non-linear approximation''), as usually
performed by compression standards (for instance, JPEG and JPEG2000).

\section{Recollection on variational compression}\label{hypf-sec}
Here and all along the paper $N$ is a non-negative integer, $I=\{1,\ldots,
N\}$ and $\BB = (\psi_i)_{i\in I}$ is a basis of $\RR^N$. We will also denote
for $\tau > 0$ (all along the paper $\tau$ denotes a non-negative real number)
and for all $k\in \ZZ$, $\tau_k = \tau (k-\frac{1}{2})$.

For any $(k_i)_{i \in I} \in \ZZ^N$, we denote
\begin{equation}\label{constraint}
\CC\left( (k_i)_{i \in I} \right) = \left\{\sum_{i\in I} u_i \psi_i, \forall
 i\in I, \tau_{k_i}\leq u_i \leq  \tau_{k_i+1}\right\}.
\end{equation}
We then consider the optimization problem
\[(\tilde P)\left( (k_i)_{i \in I} \right) : \left\{\begin{array}{l}
\mbox{minimize } f(v) \\
\mbox{under the constraint } v\in \CC\left( (k_i)_{i \in I} \right),
\end{array}\right.
\]
where $f$ is a norm, is continuously differentiable away from $0$ and its
level sets are strictly convex. In order to state Theorem \ref{LeTHM}, we also
need $f$ to be {\em curved}. This means that the inverse of the
homeomorphism\footnote{We prove in \cite{geomoptim} that, under the above
  hypotheses, $h$ actually is an homeomorphism.} $h$ below is  Lipschitz.
\[\begin{array}{rcl}
h:\{u\in\RR^N,\  f(u)=1\} & \rightarrow & \{g\in\RR^N, \
\|g\|_2=1\} \\
u & \mapsto & \frac{\nabla f(u)}{\|\nabla f(u) \|_2}.
\end{array}
\]
(The notation $\|.\|_2$ refers to the euclidean norm in $\RR^N$.)

We denote, for any $(k_i)_{i \in I} \in \ZZ^N$,
\[\tilde J \left( (k_i)_{i \in I} \right) = \{i \in I,  u^*_i =
\tau_{k_i} \mbox{ or }u^*_i = \tau_{k_i+1}  \},
\]
where $u^* = \sum_{i\in I} u^*_i \psi_i $ is the solution to
$(\tilde P)\left( (k_i)_{i \in I} \right)$. 

The interest for these optimization problems comes from the fact that, as
explained in \cite{MalgouyresCompression}, we can recover $(k_i)_{i \in I}$
from the knowledge of $(\tilde J , (u^*_i)_{j\in \tilde J} )$ (where
$\tilde J = \tilde J \left( (k_i)_{i \in I} \right)$). 

The problem $(P)$ can therefore be used for compression. Given a datum
$u=\sum_{i\in I} u_i \psi_i\in\RR^N$, we consider the unique $ (k_i(u))_{i \in
  I}\in\ZZ^N$ such that (for instance) 
\begin{equation}\label{approx}
\forall i\in I, \tau_{k_i(u)} \leq u_i <\tau_{k_i(u)+1}.
\end{equation}
The information $(\tilde J , (u^*_i)_{j\in \tilde J} )$, where $\tilde J =
\tilde J \left( (k_i(u))_{i \in I} \right)$,  is then used to encode
$u$. In the remainder, we denote  the set of indexes that need to be coded to
describe $u$ by $\tilde J(u) = \tilde J\left((k_i(u))_{i \in I}\right)$.

Notice that we can also show (see \cite{MalgouyresCompression}) that the
coding performed by the standard image processing compression algorithms (JPEG
and JPEG2000) corresponds to the above model when, for instance,
\[f(\sum_{i\in I} u_i \psi_i) = \sum_{i\in I} |u_i|^2.
\]

\section{The estimate}\label{estim-sec}
The theorem takes the form :
\begin{theorem}\label{LeTHM}
Let $\tau'>0$ and $U$ be a random variable whose low is uniform in
$\LS{\tau'}$, for a norm $f_d$. Assume $f$ satisfies the hypotheses given in
Section \ref{hypf-sec}. For any norm $\|.\|$ and any $K\in\{1,\ldots N\}$ there
exists $D_K$ such that for all $\varepsilon >0$, there exists $T>0$ such that
for all $\tau<T$
\[\PROBA{\# \tilde J\left(U \right) = K} \leq 
D_K E^{\frac{N-K}{N+1}} + \varepsilon,
\]
where $E$ is the approximation error\footnote{When computing the approximation
  error, we consider the center of $\CC\left((k_i)_{i\in I}\right)$ has
  been chosen to represent all the elements such that $
    (k_i)_{i \in I}  = (k_i(u))_{i \in I} $. } :
\[ E = \EXPECT{\|U-\tau \sum_{i\in I} k_i(U) \psi_i \| }.
\]

Moreover, if $f(\sum_{i\in I} u_i \psi_i) = \sum_{i\in I } |u_i|^2$, we also
have\footnote{This assumption is very pesimistic. For instance, the lower bound
  seems to hold for almost every basis $\BB$ of $\RR^N$, when $f$ is fixed. We
  have not worked the details of the proof of such a statement out though.}
\[\PROBA{\# \tilde J\left( U \right) =K} \geq 
D_K E^{\frac{N-K}{N+1}} - \varepsilon.
\]
\end{theorem}

When compared to the kind of results evoked in Section \ref{intro-sec}, the
above theorem original is in several ways : 

First, it concerns variational models which are more general than the model in
which the results of Section \ref{intro-sec} are usually stated. This is
probably the main interest of the current result. For instance, by following
reasonings similar to those which led to Theorem \ref{LeTHM}, it is probably
possible to obtain approximation results with redundant transforms.

Secondly, it expresses the distribution of the number of coefficients as a
function of the approximation error, while former results do the
opposite. Typically, they bound the approximation error (quantified by the
$L^2$ norm) by a function of the number of coefficients that are coded. The
comparative advantages and drawbacks of the two kind of statements is not very
clear. In the framework of Theorem \ref{LeTHM}, the larger $D_K$ (for $K$
small), the better the model compresses the data. However, it is clear that,
as the approximation error goes to $0$, we have more and more chances to
obtain a code of size $N$. With this regard, the constant $D_{K-1}$ seems to
play a particular role since it dominates (asymptotically as $\tau$ goes to
$0$) the probability not to obtain a code of length $N$.

Thirdly, it is stated in finite dimension and, as a consequence, it does not
impose apriori links between the data distribution (the function $f_d$) and
the model (the function $f$ and the basis $\BB$). The ability of the model to
represent the data is always defined. For instance, this allows the
comparison of two bad models (which is not possible in infinite dimension). 
The analog of Theorem \ref{LeTHM} in infinite dimension might be interesting,
though.

\newpage

\section{Proof of Theorem \ref{LeTHM}}

\subsection{First properties and recollection}
\subsubsection{Rewriting $(\tilde P)$}
For any $u\in \RR^N$, $(P)(u)$ denotes the optimization problem
\[(P)\left( u \right) : \left\{\begin{array}{l}
\mbox{minimize } f(v-u) \\
\mbox{under the constraint } v\in \CC\left( 0\right),
\end{array}\right.
\]
where $0$ denotes the origin in $\ZZ^N$ and $\CC(.)$ is defined by
\eqref{constraint}.

We then denote, for any $u = \sum_{i\in I}u_i \psi_i \in \CC\left(
  0\right) $,
\[J(u) = \{i \in I, u_i = \frac{\tau}{2} \mbox{ or }   u_i = - \frac{\tau}{2}
\}.
\]
With this notation, the set of active constraints of the solution $u^*$ to
$(P)\left( u \right)$ is simply $J(u^*)$.
\begin{e-proposition}
For any $(k_i)_{i \in I} \in \ZZ^N$
\[\tilde J \left( (k_i)_{i \in I} \right) = J(u^*), 
\]
where $u^*$ is the solution to $(P)\left(\tau \sum_{i \in I}k_i \psi_i
\right)$.
\end{e-proposition}
\proof
Denoting $\tilde u^*$ the solution of $(\tilde P)\left( (k_i)_{i \in I}
\right)$ and $u^*$ the solution to $(P)\left( \tau \sum_{i \in I}k_i \psi_i
\right)$, we have
\begin{equation}\label{equivP}
\tilde u^* = u^* + \tau \sum_{i \in I}k_i \psi_i.
\end{equation}
This can be seen from the fact that $(P)\left( \sum_{i \in I}k_i \psi_i
\right)$ is exactly $(\tilde P)\left( (k_i)_{i \in I}
\right)$, modulo a "global translation" by $\tau \sum_{i \in I}k_i
\psi_i$. (The rigorous proof of \eqref{equivP} can easily be established
using Kuhn-Tucker conditions, see \cite{Rockafellar}, Th 28.3, pp. 281.)

The proposition is then obtained by identifying the coordinates of $\tilde
u^*$ and $u^*$ in the basis $\BB$.
\endproof
\subsubsection{On projection onto polytopes}\label{proj-sec}
We can now adapt the definitions and notations of \cite{geomoptim} to the
problems $(P)(.)$. Beside Proposition \ref{ssat-prop1}, all the results stated
in this section are proved in \cite{geomoptim}.

We consider a norm $f_d$ (which will be used latter on to define the data
distribution law) and define for any $C\subset \RR^N$  and any $A\subset \RR$
\[\SSAT{C}{A} = \left\{u \in \RR^N, \exists u^* \in C, u^*\mbox{ is solution
    to } (P)(u) \mbox{ and } f_d(u-u^*) \in A \right\}.
\]
This corresponds to all the optimization problems whose solution is in $C$ (we
also control the distance between $u$ and the result of $(P)(u)$). Notice that
$\SSAT{C}{A}$ depends on $\tau$. We do not make this dependence explicit since
it does not create any confusion, in practice.

We also define the equivalence relationship over $\CC(0)$
\[u\sim v \Longleftrightarrow J(u) = J(v).
\]
For any $u\in \CC(0)$, we denote $\overline u$ the equivalence class of $u$.

In the context of this paper, we obviously have for all $u =  \sum_{i\in
  I}u_i \psi_i \in \CC(0)$ 
\begin{equation}\label{ubarre}
\overline u = \left\{u^c + \tau \sum_{j\not\in J(u)}\beta_j \psi_j, \forall
  j\not\in J(u), -\frac{1}{2} <\beta_j <\frac{1}{2} \right\},
\end{equation}
where
\[u^c = \sum_{j\in J(u)} u_j \psi_j.
\]
(Here and all along the paper the notation $j\not\in J$ stands for $j\in
I\setminus J$.)

Let us give some descriptions of $\SSAT{.}{.}$.

\begin{e-proposition}\label{ssat-transl}
For any $u^*\in \partial \CC(0)$ and any $v\in \overline{u^*}$,
\[\SSAT{v}{1} = (v-u^*) + \SSAT{u^*}{1}.
\]
\end{e-proposition}
In words, $\SSAT{v}{1}$ is a translation of $\SSAT{u^*}{1}$.

\begin{e-proposition}\label{ssat-rescale}
For any $u^*\in \partial \CC(0)$, any $v\in \SSAT{u^*}{]0,+\infty[}$
    and any $\lambda>0$
\[u^*+\lambda (v-u^*)  \in \SSAT{u^*}{]0,+\infty[}.
\]
\end{e-proposition}

\begin{theorem}\label{thm1}
For any $u^* \in \partial \CC(0)$ and any $\tau' >0$,
\[\SSAT{\overline{u^*}}{]0,\tau']} = \left\{v + \lambda(u-u^*), \mbox{for }
  v\in \overline{u^*}, \lambda \in ]0,\tau'] \mbox{ and } u \in \SSAT{u^*}{1}
\right\}
\]
\end{theorem}

%

We also have (see \cite{geomoptim})
\begin{e-proposition}\label{ssat-prop0}
If $f$ satisfies the hypotheses given in Section \ref{hypf-sec},
for any $u^*\in\partial \CC(0)$, $\SSAT{u^*}{1}$ is a non-empty, compact
Lipschitz manifold of dimension $\#J(u^*)-1$.
\end{e-proposition}

Another useful result for the purpose of this paper is the following.
\begin{e-proposition}\label{ssat-prop1}
If $f$ satisfies the hypotheses given in Section \ref{hypf-sec},
for any $u^*\in\partial \CC(0)$ and any $\tau'>0$,
$\SSAT{u^*}{]0,\tau']}$ is a non-empty, bounded Lipschitz manifold of
dimension $\#J(u^*)$.
\end{e-proposition}

\proof
In order to prove the proposition, we consider $u^*=\sum_{i\in I} u^*_i \psi_i
\in \partial\CC(0)$ and $u^c =  \sum_{i\in J(u^*)} u^*_i \psi_i$. We are
going to prove the proposition in the particular case where
$u^c=u^*$. Proposition \ref{ssat-transl} and \ref{ssat-rescale} permit indeed
to generalize the latter result obtained to any $\SSAT{u^*}{]0,\tau']}$, for
$u^*\in \overline{u^c}$. (They indeed guarantee that $\SSAT{u^*}{]0,\tau']}$ is
obtained by translating $\SSAT{u^c}{]0,\tau']}$.)

In order to prove that $\SSAT{u^c}{]0,\tau']}$ is a bounded Lipschitz manifold
of dimension $\#J(u^*)$, we prove that the mapping $h'$ defined below is a
Lipschitz homeomorphism.
\begin{equation}\label{hprime}
\begin{array}{rrcl}
h' : & \SSAT{u^c}{1}\times ]0,\tau'] & \longrightarrow &
    \SSAT{u^c}{]0,\tau']}  \\
     &    (u,\lambda) & \longmapsto & u^c + \lambda (u-u^c). 
\end{array}
\end{equation}
The conclusion then directly follows from Proposition \ref{ssat-prop0}.

Notice first that we can deduce from Proposition \ref{ssat-rescale}, that $h'$
is properly defined.

Let us prove that $h'$ is invertible. For this purpose, we consider
$\lambda_1$ and $\lambda_2$ in $]0,\tau']$ and $u_1$ and $u_2$ in
    $\SSAT{u^c}{1}$ such that
\begin{equation}\label{tmp98656}
u^c+\lambda_1 (u_1-u^c) = u^c+\lambda_2 (u_2-u^c).
\end{equation}
We have
\begin{eqnarray*}
\lambda_1 & = & f_d(\lambda_1 (u_1-u^c)) \\
 & = & f_d(\lambda_2 (u_2-u^c)) \\
& = & \lambda_2. 
\end{eqnarray*}
Using \eqref{tmp98656}, we also obtain $u_1=u_2$ and $h'$ is
invertible.

Finally, $h'$ is Lipschitz since, for any $\lambda_1$ and $\lambda_2$
    in $]0,\tau']$ and any $u_1$ and $u_2$ in $ \SSAT{u^c}{1}$,
\begin{eqnarray*}
\|\lambda_1 (u_1-u^c) - \lambda_2 (u_2-u^c) \|_2 & = &  \|\lambda_1
(u_1-u_2) + (\lambda_1-\lambda_2) (u_2-u^c)  \|_2, \\
& \leq & \tau' \|u_1-u_2\|_2 + C |\lambda_1-\lambda_2|,
\end{eqnarray*}
where $C$ is such that for all $u\in \SSAT{u^c}{1}$,
\[\|u-u^c\|_2 \leq C.
\]
(Remember $\SSAT{u^c}{1}$ is compact, see Proposition \ref{ssat-prop0}.)
\endproof

\subsection{The estimate}
We denote the discrete grid by
\[\DD = \{\tau \sum_{i\in I} k_i \psi_i, (k_i)_{i\in I} \in \ZZ^N\},
\]
and, for $u^*\in\partial \CC(0)$ and $(k_j)_{j\in J(u^*)}\in
\ZZ^{J(u^*)}$,
\[\Dom{ (k_j)_{j\in J(u^*)} }= \{\tau \sum_{j\in J(u^*)} k_j \psi_j
+\tau \sum_{i\not\in J(u^*)} k_i \psi_i  , \mbox{ where } (k_i)_{i\not\in
  J(u^*)} \in \ZZ^{I\setminus J(u^*)}\}.
\]
The set $\Dom{ (k_j)_{j\in J(u^*)} }$ is a slice in $\DD$.
\begin{e-proposition}\label{major1}
Let $\tau'>0$, $u^*\in\partial \CC(0)$ and $(k_j)_{j\in J(u^*)}\in
\ZZ^{J(u^*)}$,
\[\#\left(   \SSAT{\overline{u^*}}{]0,\tau']}  \cap \Dom{ (k_j)_{j\in
	J(u^*)} }\right) \leq 1.
\]
\end{e-proposition}
\proof
Taking the notations of the proposition and assuming $
\SSAT{\overline{u^*}}{]0,\tau']}  \cap \Dom{ (k_j)_{j\in J(u^*)}} \neq
      \emptyset$, we consider $(k_i^1)_{i\in I}$ and  $(k_i^2)_{i\in I}$ such
      that
\[\tau  \sum_{i\in I} k_i^1 \psi_i \in
\SSAT{\overline{u^*}}{]0,\tau']}  \cap \Dom{ (k_j)_{j\in J(u^*)} }
\] 
and
\[\tau  \sum_{i\in I} k_i^2 \psi_i \in
\SSAT{\overline{u^*}}{]0,\tau']}  \cap \Dom{ (k_j)_{j\in J(u^*)} }.
\]
Theorem \ref{thm1} guarantees there exist $v_1$ and $v_2$ in
$\overline{u^*}$, $\lambda_1$ and $\lambda_2$ in $]0,\tau']$ and $u_1$
    and $u_2$ in $\SSAT{u^*}{1}$ such that
\[\tau \sum_{i\in I} k_i^1 \psi_i = v_1+\lambda_1(u_1-u^* )
\] 
and
\[\tau \sum_{i\in I} k_i^2 \psi_i = v_2+\lambda_2(u_2-u^* ).
\] 
So
\[v_1+\lambda_1(u_1-u^* ) = v_2+\lambda_2(u_2-u^* ) + \tau
\sum_{i\not\in J(u^*)} ( k_i^1 - k_i^2 )\psi_i.
\]

Using \eqref{ubarre}, we know there exists $(\beta_i^1)_{i\not\in
  J(u^*)}$ and $(\beta_i^2)_{i\not\in J(u^*)}$ such that
\[\forall i\not\in J(u^*), -\frac{1}{2}< \beta_i^1 <\frac{1}{2}
  \mbox{ and } -\frac{1}{2}< \beta_i^2 <\frac{1}{2},
\] 
\[v_1 = u^c + \tau \sum_{i\not\in J(u^*)} \beta_i^1 \psi_i
\]
and
\[v_2 = u^c + \tau \sum_{i\not\in J(u^*)} \beta_i^2 \psi_i,
\]
with $u^c = \sum_{j\in J(u^*)} u^*_j \psi_j$, where $u^* = \sum_{i\in I} u^*_i
\psi_i$.

So, letting for all $i\not\in J(u^*)$, $\alpha_i =k_i^1 - k_i^2  +
\beta_i^2 -  \beta_i^1$, we finally have
\begin{equation}\label{tmp678}
\lambda_1(u_1-u^* ) = \lambda_2(u_2-u^* ) + \tau
\sum_{i\not\in J(u^*)} \alpha_i \psi_i.
\end{equation}
Let us assume 
\begin{equation}\label{hyp}
\max_{i\not\in J(u^*)}  |\alpha_i| >0,
\end{equation}
 and consider $0<\lambda \leq 1$ such that
\begin{equation}\label{lambda}
\lambda < \frac{1}{2 \max_{i\not\in J(u^*)}  |\alpha_i| }.
\end{equation}
We have, using \eqref{tmp678},
\begin{eqnarray*}
u^c + \lambda\lambda_1 [(u_1-u^*+u^c)-u^c ] & = & u^c + \lambda\lambda_1
(u_1-u^* ) \\
& = & u^c + \lambda \tau \sum_{i\not\in J(u^*)} \alpha_i \psi_i + \lambda
\lambda_2(u_2-u^* ) \\
& = & v + \lambda \lambda_2  [(u_2-u^*+v)-v ],
\end{eqnarray*}
where $v= u^c + \lambda\tau \sum_{i\not\in J(u^*)} \alpha_i \psi_i$. Moreover,
using \eqref{ubarre} and \eqref{lambda}, we know that $v\in
\overline{u^c}$. Using Proposition \ref{ssat-transl}, we know that
\[u_1-u^*+u^c \in \SSAT{u^c}{1} \mbox{ and } u_2-u^*+v \in \SSAT{v}{1}.
\]
Finally, applying Theorem \ref{thm1}, we obtain 
\[u^c + \lambda\lambda_1 (u_1-u^* ) \in \SSAT{u^c}{]0,\tau']}\cap
    \SSAT{v}{]0,\tau']}.
\]
Since the solution to $(P)(u^c + \lambda\lambda_1 (u_1-u^* ))$ is
unique, we necessarily have $u^c=v$ and therefore $\max_{i\not\in J(u^*)}
|\alpha_i| =0$. This contradicts \eqref{hyp} and guarantees that 
\[\max_{i\not\in J(u^*)}  |\alpha_i| =0.
\]
Using the definition of $\alpha_i$, we obtain, for all $i\not\in J(u^*)$,
\[|k_i^1 - k_i^2|  = |\beta_i^1 -  \beta_i^2| <1.
\]
This implies $k_i^1 = k_i^2$, for all $i\in I$.
\endproof
Let us denote, for $u^*\in\partial \CC(0)$, the projection onto $Span\left(
  \psi_j, j\in J(u^*)\right)$ by
\[\begin{array}{rrcl}
p : & \RR^N & \longrightarrow & Span\left( \psi_j, j\in J(u^*)\right)\\
     & \sum_{i\in I} \alpha_i \psi_i & \longmapsto & \sum_{j\in J(u^*) }
     \alpha_j \psi_j.
\end{array}
\]
It is not difficult to see that, for any $\tau'>0$, $u^*\in\partial \CC(0)$ and
$(k_j)_{j\in J(u^*)}\in \ZZ^{J(u^*)}$,
\begin{equation}\label{bbbb}
\#\left( \SSAT{\overline{u^*}}{]0,\tau']}  \cap \Dom{ (k_j)_{j\in J(u^*)}
    }\right) = 1 \Longrightarrow \tau \sum_{j\in J(u^*)} k_j\psi_j\in p\left(
\SSAT{\overline{u^*}}{]0,\tau']}  \right).
\end{equation}
\begin{remark} Notice that the converse implication does not hold in
  general. It is indeed possible to build counter examples where
$\SSAT{\overline{u^*}}{]0,\tau']} $ passes between the points of the discrete
grid $\DD$. However, it is not difficult to see that, if $\tau \sum_{j\in
  J(u^*)} k_j\psi_j\in p\left( \SSAT{\overline{u^*}}{]0,\tau']}  \right)$ and 
$\SSAT{\overline{u^*}}{]0,\tau']}  \cap \Dom{ (k_j)_{j\in J(u^*)} }=
\emptyset$, we can build $(k_i)_{i\not\in J(u^*)} \in \ZZ^{J\setminus J(u^*)}$
such that
\[\tau \sum_{j\in J(u^*)} k_j\psi_j +\tau \sum_{i\not\in J(u^*)}
(k_i+\frac{1}{2}) \psi_i \in \SSAT{u^c}{]0,\tau']}, 
\]
where
\[u^c = \sum_{j\in J(u^*)} u^*_j \psi_j. (u^*_j \mbox{ are the coordinates of
}u^*)
\]
This means that the set $\SSAT{u^c}{]0,\tau']}$, which is a manifold of
dimension $\# J(u^c)$ living in $\RR^N$, intersects a discrete grid. This is
obviously a very rare event. Typically, adding to the basis $\BB$ some kind
of randomness (for instance adding a very small Gaussian noise to every
$\psi_i$) would make it an event of probability $0$.

Notice, with this regard, that when $f(\sum_{i\in I} u_i \psi_i) = \sum_{i \in
I} |u_i|^2$, we trivially have the equivalence in \eqref{bbbb}.
\end{remark}

A simple consequence of \eqref{bbbb} is that
\begin{equation}\label{borne-1}
\#\left( \SSAT{\overline{u^*}}{]0,\tau']}  \cap \DD\right) \leq
    \#\left(p\left( \SSAT{\overline{u^*}}{]0,\tau']} \right)  \cap \left\{\tau
    \sum_{j\in J(u^*)} k_j\psi_j, (k_j)_{j\in J(u^*)}\in
\ZZ^{J(u^*)}\right\} \right). 
\end{equation}
Notice finally that, for $u^*=\sum_{i\in I} u^*_i\psi_i \in\partial
\CC(0)$, Proposition \ref{ssat-transl} and Equation \eqref{ubarre} guarantees
that
\[p\left(\SSAT{u^c}{1}\right) = p\left(\SSAT{u^*}{1}\right), 
\]
for $u^c = \sum_{j\in J(u^*)} u^*_j \psi_j$.

We therefore have, using also Theorem \ref{thm1}, Proposition
\ref{ssat-rescale} and Equation \eqref{ubarre},
\begin{eqnarray*}
p\left( \SSAT{\overline{u^*}}{]0,\tau']} \right) & = &\{ p(v) +\lambda
    (p(u) - p(u^*)), \mbox{ for } v\in \overline{u^*}, \lambda\in ]0,\tau']
    \mbox{ and } u\in \SSAT{u^*}{1} \}, \\
 & = &\{ u^c +\lambda (p(u) - u^c),  \mbox{ for } \lambda\in ]0,\tau'] \mbox{
   and } u\in \SSAT{u^c}{1}\}, \\
& = & p\left( \SSAT{u^c}{]0,\tau']} \right).
\end{eqnarray*}

Finally,
\begin{equation}\label{borne0}
\#\left( \SSAT{\overline{u^*}}{]0,\tau']}  \cap \DD\right) \leq
    \#\left(p\left( \SSAT{u^c}{]0,\tau']} \right)  \cap \left\{\tau
    \sum_{j\in J(u^c)} k_j\psi_j, (k_j)_{j\in J(u^c)}\in
\ZZ^{J(u^c)}\right\} \right). 
\end{equation}

\begin{e-proposition}\label{mesurable}
If $f$ satisfies the hypotheses given in Section \ref{hypf-sec}
then, for any $u^* = \sum_{i\in I} u^*_i \psi_i \in\partial \CC(0)$,
$p\left(\SSAT{u^c}{]0,\tau']}\right)$ (where $u^c = \sum_{j\in J(u^*)} u^*_j
\psi_j$) is a non-empty, bounded Lipschitz manifold of dimension $\#J(u^*)$. 
\end{e-proposition}
\proof
Thanks to Proposition \ref{ssat-prop1}, it suffices to establish that the
restriction of $p$ :
\[\begin{array}{rrcl}
p' : & \SSAT{u^c}{]0,\tau']} & \longrightarrow &
p\left(\SSAT{u^c}{]0,\tau']}\right)\\
     &    u & \longmapsto & p(u).
\end{array}
\]
is a Lipschitz homeomorphism. This latter result is immediate once
we have established that $p'$ is invertible. 

This proof is similar to the one of Proposition \ref{major1}. Taking
the notations of the proposition, we assume
that there exist $u_1$ and $u_2$ in $ \SSAT{u^c}{]0,\tau']}$ and
$(\alpha_i)_{i\not\in J(u^*)}\in \RR^{J(u^*)}$  satisfying
\[u_1 = u_2 + \tau \sum_{i\not\in J(u^*)} \alpha_i \psi_i.
\]
If we assume $\max_{i\not\in J(u^*)} |\alpha_i|\neq 0$, we have for $0<\lambda
< \min(1,\frac{1}{2 \max_{i\not\in J(u^*)} |\alpha_i|})$,
\begin{eqnarray*}
u^c+\lambda(u_1-u^c) & = & u^c + \tau \sum_{i\not\in J(u^*)} \lambda
\alpha_i \psi_i + \lambda (u_2-u^c) \\
& = & v + \lambda \left(u_2 + \tau \sum_{i\not\in J(u^*)} \lambda
\alpha_i \psi_i -v\right)
\end{eqnarray*}
for $v =  u^c + \tau \sum_{i\not\in J(u^*)} \lambda \alpha_i \psi_i
$. Since $v\in \overline{u^c}$ (see \eqref{ubarre}), Proposition
\ref{ssat-transl} guarantees that $u_2 + \tau \sum_{i\not\in J(u^*)} \lambda
\alpha_i \psi_i= u_2+v-u^c\in \SSAT{v}{]0,\tau]} $. As a consequence, applying
Proposition \ref{ssat-rescale}, we know that
\[u^c+\lambda (u_1-u^c) \in \SSAT{u^c}{\lambda} \cap  \SSAT{v}{]0,+\infty[}.
\]

Since $(P)(u^c+\lambda (u_1-u^c))$ has a unique solution, we obtain a
contradiction and can conclude that for all $i\not\in J(u^*)$, $\max_{i\not\in
  J(u^*)} |\alpha_i|= 0$. 

As a consequence, $p'$ is invertible. It is then obviously a Lipschitz
homeomorphism.
\endproof

Proposition \ref{mesurable} guarantees that
$p\left(\SSAT{u^c}{]0,\tau']}\right)$ is Lebesgue measurable in
$\RR^{\#J(u^*)}$. Moreover, its Lebesgue measure in $\RR^{\#J(u^*)}$ (denoted
$\LEB{\#J(u^*)}{p\left(\SSAT{u^c}{]0,\tau']}\right)}$) is finite and strictly
positive : 
\[0< \LEB{\#J(u^*)}{p\left(\SSAT{u^c}{]0,\tau']}\right)} <\infty.
\]

Another consequence takes the form of the following proposition.
\begin{e-proposition}\label{estime1}
Let $\tau'>0$ and $u^*\in\partial \CC(0)$
\[\lim_{\tau\rightarrow 0} \tau^K \#\left( \SSAT{\overline{u^*}}{]0,\tau']}
  \cap \DD\right) \leq \LEB{K}{p\left(\SSAT{u^c}{]0,\tau']}\right)}
\]
where $K=\#J(u^*)$.

Moreover, if the equality holds in \eqref{borne-1} (or equivalently : the
equality holds in \eqref{borne0}) 
\[\lim_{\tau\rightarrow 0} \tau^K \#\left( \SSAT{\overline{u^*}}{]0,\tau']}
  \cap \DD\right) = \LEB{K}{p\left(\SSAT{u^c}{]0,\tau']}\right)}.
\]
\end{e-proposition}
\proof
In order to prove the proposition, we are going to prove that, denoting
$K=\#J(u^c)$,
\begin{equation}\label{09u82t}
\lim_{\tau\rightarrow 0} \tau^{K} \#\left(p\left(
    \SSAT{u^c}{]0,\tau']} \right)  \cap \left\{\tau \sum_{j\in J(u^c)}
    k_j\psi_j, (k_j)_{j\in J(u^c)}\in \ZZ^{J(u^c)}\right\} \right) =
\LEB{K}{p\left(\SSAT{u^c}{]0,\tau']}\right)} 
\end{equation}
The conclusion follows from \eqref{borne0}.

Let us first remark that, unlike $\SSAT{u^c}{]0,\tau']}$,  the set
\[A = p\left(\SSAT{u^c}{]0,\tau']}\right) - u^c
\]
does not depend on $\tau$. This is due to Proposition 9\footnote{The
  definition of $\SSAT{C}{A}$ given in the current paper does not allow the
  rewriting of the proposition 9 of \cite{geomoptim}. This is why we have not
  adapted it in Section \ref{proj-sec}.}, in \cite{geomoptim}. Notice also
that, because of Proposition \ref{mesurable}, both $A$ and
$p\left(\SSAT{u^c}{]0,\tau']}\right)$ are Lebesgue measurable (in $\RR^K$) and
that
\[\LEB{K}{A} = \LEB{K}{p\left(\SSAT{u^c}{]0,\tau']}\right)}. 
\]

In order to prove the upper bound in \eqref{09u82t}, we consider the sequence
of functions, defined over $\RR^K$
\[f_n(u) = \max\left(0 , 1-n \ \inf_{v\in A} \|u-v\|_2 \right).
\]

This is a sequence of functions which are both Lebesgue and Riemann integrable
and the sequence converges in $L^1(\RR^K)$ to $\one_{A}$ (the indicator
function of the set $A$). So, for any $\varepsilon>0$, there exists $n\in \NN$
such that
\[\int f_n \leq \int \one_A + \varepsilon.
\]
 Moreover, we have, for all $u\in\RR^K$ and all $n\in\NN$, 
\[\one_{A}(u) \leq f_n(u).
\]
So, denoting $V_{\tau}= \left\{\tau \sum_{j\in J(u^c)} k_j\psi_j - u^c,
  (k_j)_{j\in J(u^c)}\in \ZZ^{J(u^c)}\right\}$,
\begin{eqnarray*}
\lim_{\tau\rightarrow 0}\tau^{K} \#\left(p\left(
    \SSAT{u^c}{]0,\tau']} \right)  \cap \left\{\tau \sum_{j\in J(u^c)}
    k_j\psi_j, (k_j)_{j\in J(u^c)}\in \ZZ^{J(u^c)}\right\} \right)   & =  & 
\lim_{\tau\rightarrow 0}\tau^{K} \sum_{v\in V_{\tau} }\one_{A}(v) \\
& \leq  & \lim_{\tau\rightarrow 0}\tau^{K} \sum_{v\in V_{\tau} }  f_n(v) \\
& \leq  & \int f_n \\
& \leq  & \int \one_A + \varepsilon \\
& \leq  &  \LEB{K}{p\left(\SSAT{u^c}{]0,\tau']}\right)} + \varepsilon.
\end{eqnarray*}
So,
\[\lim_{\tau\rightarrow 0}\tau^{K} \#\left(p\left(
    \SSAT{u^c}{]0,\tau']} \right)  \cap \left\{\tau \sum_{j\in J(u^c)}
    k_j\psi_j, (k_j)_{j\in J(u^c)}\in \ZZ^{J(u^c)}\right\} \right) 
\leq \LEB{K}{p\left(\SSAT{u^c}{]0,\tau']}\right)}
\]

The lower bound in \eqref{09u82t} is obtained in a similar way, by considering
an approximation of $\one_{A}$ by a function smaller than $\one_{A}$ which is
Riemann integrable. (For instance : $f_n(u) = 1 - \max\left(0 , 1-n \ 
  \inf_{v\not\in A} \|u-v\|_2 \right))$.)
\endproof

From now on , we will denote for all $K\in\{1,\ldots, N\}$
\[C_K = \left\{\tau\sum_{j\in J} u_j \psi_j, \mbox{ where } J\subset I, \#J=K
  \mbox{ and } \forall j\in J, u_j =-\frac{1}{2} \mbox{ or } u_j =\frac{1}{2}
\right\}
\]
The set $C_K$ contains all the "centers" of the equivalence classes of
codimension $K$.

Similarly, we denote
\[\Cl{K} = \left\{u^* \in \partial \CC(0), \#J(u^*) = K \right\}.
\]
We obviously have, for all  $K\in\{1,\ldots, N\}$,
\[\Cl{K} = \cup_{u^c\in C_K} \overline{u^c}.
\]

Since, for all $K\in\{1,\ldots, N\}$, $C_K$ is finite, it is clear from
Proposition \ref{estime1} that, for any $\tau'>0$,
\[\lim_{\tau\rightarrow 0} \tau^K \#\left( \SSAT{\Cl{K}}{]0,\tau']}
  \cap \DD\right) \leq \sum_{u^c \in C_K}
\LEB{K}{p\left(\SSAT{u^c}{]0,\tau']}\right)} < +\infty
\]
Moreover, we have an equality between the above two terms, as soon as the
equality holds in \eqref{borne-1}. 

We can finally express the following estimate.
\begin{e-proposition}\label{finalprop}
Let $\tau'>0$
\[\lim_{\tau\rightarrow 0} \tau^K \#\left( \SSAT{\Cl{K}}{]0,\infty[}
  \cap \LS{\tau'}\cap \DD\right) \leq \sum_{u^c \in C_K}
\LEB{K}{p\left(\SSAT{u^c}{]0,\tau']}\right)}
\]
where $K=\#J(u^*)$.

Moreover, if the equality holds in \eqref{borne-1} for all $u^c\in C_K$ 
(or equivalently : the equality holds in \eqref{borne0}) 
\[\lim_{\tau\rightarrow 0} \tau^K \#\left( \SSAT{\Cl{K}}{]0,\infty[}\cap
  \LS{\tau'} \cap \DD\right) =  \sum_{u^c \in
  C_K}\LEB{K}{p\left(\SSAT{u^c}{]0,\tau']}\right)}. 
\]
\end{e-proposition}
\proof
We consider 
\[M = \sup_{\{u=\sum_{i\in I} u_i \psi_i, \forall i\in I, |u_i|\leq
  \frac{1}{2}\}} f_d(u)
\]
We have, for all $u^*\in\partial \CC(0)$,
\begin{equation}\label{M}
f_d(u^*) \leq M \tau.
\end{equation}

We therefore have for all $u\in \LS{\tau'}$ and for $u^*$ the solution to
$(P)(u)$,
\begin{eqnarray*}
f_d(u-u^*) & \leq & f_d(u) + f_d(u^*) \\
& \leq & \tau' + M\tau.
\end{eqnarray*}
So
\[\SSAT{\Cl{K}}{]0,\infty[}\cap \LS{\tau'} \subset \SSAT{\Cl{K}}{]0,\tau'+M
  \tau]}.
\]
Moreover, it is not difficult to see that (remember $h'$ defined by
\eqref{hprime} is an homeomorphism)
\[\lim_{\tau\rightarrow 0} \sum_{u^c \in C_K} \LEB{K}{p\left(
    \SSAT{u^c}{]0,\tau'+M \tau]}\right)} =  \sum_{u^c \in C_K} \LEB{K}{p\left(
    \SSAT{u^c}{]0,\tau']}\right)}.
\]

We can therefore deduce (from Proposition \ref{estime1}) that 
\[\lim_{\tau\rightarrow 0} \tau^K \#\left( \SSAT{\Cl{K}}{]0,\infty[}
  \cup \LS{\tau'}\cap \DD\right) \leq \sum_{u^c \in C_K}
\LEB{K}{p\left(\SSAT{u^c}{]0,\tau']}\right)}
\]

In order to prove the last statement of the proposition, we consider $u^*\in
\partial \CC(0)$ and $u\in\SSAT{u^*}{]0,\tau']}$, we know that
\begin{eqnarray*}
f_d(u) & \leq &  f_d(u-u^*) + f_d(u^*) \\
  & \leq &  \tau'+ M\tau \\
\end{eqnarray*}

So
\[ \SSAT{\Cl{K}}{]0,\tau'-M \tau]} \subset\SSAT{\Cl{K}}{]0,\infty[}\cap
\LS{\tau'}.
\]
Since (again)
\[\lim_{\tau\rightarrow 0} \sum_{u^c \in C_K} \LEB{K}{p\left(
    \SSAT{u^c}{]0,\tau'-M \tau]}\right)} =  \sum_{u^c \in C_K} \LEB{K}{p\left(
    \SSAT{u^c}{]0,\tau']}\right)},
\]
we know that the second statement of the proposition holds.
\endproof

Another immediate result is useful to state the final theorem. Notice first
that we have, for any $(k_i)_{i\in I} \in \ZZ^N$ and any norm $\|.\|$,
\[\int_{v\in \CC((k_i)_{i\in I})} \|v-\tau \sum_{i\in I} k_i \psi_i \| dv = C
\tau^{N+1},
\]
where 
\[C = \int_{\{v = \sum_{i\in I} v_i \psi_i, \forall i\in I, |v_i|\leq
  \frac{1}{2}\}} \|v\| dv
\] 
only depends on the particular norm $\|.\|$ and the basis
$(\psi_i)_{i\in I}$.

So, denoting $U$ a random variable whose law is uniform in $\LS{\tau'}$ and
$(k_i(U))_{i\in I}$ the discrete point defined by \eqref{approx}, we have
\begin{equation}\label{errorespect}
\lim_{\tau\rightarrow 0} \frac{\EXPECT{\|U-\tau \sum_{i\in I} k_i(U) \psi_i \|
  } } {\tau^{N+1}} = C.
\end{equation}

This follows from the fact that the number of points $(k_i)_{i\in I}$ such
that $\CC((k_i)_{i\in I})$ intersects both $\LS{\tau'}$ and its complement in
$\RR^N$ becomes negligible with regard to the number of points $(k_i)_{i\in
  I}$ such that $\CC((k_i)_{i\in I})$ is included in $\LS{\tau'}$, when $\tau$
goes to $0$.

We can now state the final result.
\begin{theorem}
Let $\tau'>0$ and $U$ be a random variable whose low is uniform in
$\LS{\tau'}$, for a norm $f_d$. For any norm $\|.\|$, any $K\in\{1,\ldots N\}$
and any $\varepsilon >0$, there exists $T>0$ such that for all $\tau<T$
\[\PROBA{\# \tilde J\left(U \right) =K} \leq 
D_K E^{\frac{N-K}{N+1}} + \varepsilon,
\]
where $E$ is the approximation error\footnote{When computing the approximation
  error, we consider the center of $\CC\left((k_i)_{i\in I}\right)$ has
  been chosen to represent all the elements coded by $(\tilde P)\left(
    (k_i)_{i \in I} \right)$. } :
\[ E = \EXPECT{\|U-\tau \sum_{i\in I} k_i(U) \psi_i \| },
\]

Moreover, if the equality holds in \eqref{borne-1} (or
equivalently : the equality holds in \eqref{borne0}) for all $u^c\in C_K$,
then we also have
\[\PROBA{\# \tilde J\left( U \right) =K} \geq 
D_K E^{\frac{N-K}{N+1}} - \varepsilon.
\]

The constant $D_K$ is given by
\[D_K=\frac{A_K}{BC^{\frac{N-K}{N+1}}},
\]
with
\[A_K= \sum_{u^c \in C_K} \LEB{K}{p\left(\SSAT{u^c}{]0,\tau']}\right)},
\]
\[B= \frac{ \LEB{N}{\LS{\tau'}} } {\LEB{N}{\{v = \sum_{i\in I} v_i \psi_i,
    \forall i \in I, |v_i|\leq \frac{1}{2}\}} }
\]
and
\[C = \int_{\{v = \sum_{i\in I} v_i \psi_i, \forall i \in I, |v_i|\leq
  \frac{1}{2}\}} \|v\| dv.
\]
\end{theorem}
\proof
Remark first that, for any $(k_i)_{i\in I}\in \ZZ^N$, the probability that
\[\tau_{k_i} \leq U_i \leq \tau_{k_i+1},
\]
when $U=\sum_{i\in I} U_i \psi_i$ follows a uniform law in $\LS{\tau'}$, is
\[\frac{\LEB{N}{\CC((k_i)_{i\in I}) \cap \LS{\tau'} }}{\LEB{N}{\LS{\tau'}}}.
\]

Therefore, taking the notation of the theorem
\[\PROBA{\# \tilde J(U) = K } =  \sum_{(k_i)_{i\in I}\in \ZZ^N}
\one_{\tau\sum_{i\in I} k_i \psi_i \in \SSAT{\Cl{K}}{[0,+\infty[} }
\frac{\LEB{N}{\CC((k_i)_{i\in I}) \cap \LS{\tau'} }}{\LEB{N}{\LS{\tau'}}}.
\]

If $(k_i)_{i\in I}$ is such that $\LEB{N}{\CC((k_i)_{i\in I}) \cap \LS{\tau'}
} \neq 0$, there exists $v\in \CC(0)$ such that $v + \tau \sum_{i\in I} k_i
\psi_i \in \LS{\tau'}$. So, we have
\begin{eqnarray*}
f_d(\tau \sum_{i\in I} k_i \psi_i) & \leq & \tau' + f_d(v) \\
& \leq & \tau' + M\tau ,
\end{eqnarray*} 
where $M$ is given by \eqref{M}.

We therefore have
\[\PROBA{\# \tilde J(U) = K } \leq \frac{\LEB{N}{\CC(0)}}{\LEB{N}{\LS{\tau'}}}
\# \left( \SSAT{\Cl{K}}{ ]0,+\infty [ }  \cap \LS{\tau'+M\tau} \cap \DD
\right).
\]

The lower bound is obtained with a similar estimation and we obtain
\[\PROBA{\# \tilde J(U) = K } \geq \frac{\LEB{N}{\CC(0)}}{\LEB{N}{\LS{\tau'}}}
\# \left( \SSAT{\Cl{K}}{]0,+\infty[} \cap \LS{\tau'-M\tau} \cap \DD \right).
\]

Notice finally that
\[\lim_{\tau\rightarrow 0} \frac{\# \left( \SSAT{\Cl{K}}{]0,+\infty[} \cap
    \LS{\tau'} \cap \DD \right)}{\# \left( \SSAT{\Cl{K}}{]0,+\infty[} \cap
    \LS{\tau'\pm M\tau} \cap \DD \right)} =1.
\]

The proof is now a straightforward consequence of Proposition \ref{finalprop}
and \eqref{errorespect}. More precisely, taking the notations of
the theorem and  $\varepsilon>0$, we know that there exists $T>0$ such that,
for all $\tau<T$,
\[ \tau^K \#\left( \SSAT{\Cl{K}}{]0,\infty[}
  \cup \LS{\tau'+M\tau}\cap \DD\right) \leq A_K + \varepsilon,
\]
and
\[\frac{E^{\frac{1}{N+1}} } {C^{\frac{1}{N+1}}} \geq \tau - \varepsilon.
\]

So
\begin{eqnarray*}
\PROBA{\# \tilde J\left((K_i)_{i\in I} \right) =K}
 & \leq & \frac{\tau^N}{B} \frac{A_K+\varepsilon}{\tau^K} \\
& \leq &\frac{A_K+\varepsilon}{B}\left( \left(\frac{E}{C}
  \right)^{\frac{1}{N+1}} + \varepsilon \right)^{N-K} \\
& \leq & \frac{A_K}{B C^{\frac{N-K}{N+1}}} E^{\frac{N-K}{N+1}} +
o(1),
\end{eqnarray*}
where $o(1)$ is a function of $\varepsilon$ which goes to $0$, when
$\varepsilon$ goes to $0$. The first inequality of the theorem follows.

The proof of the second inequality of the theorem is similar to one above.
\endproof

\bibliographystyle{plain}
\markboth{}{}
\bibliography{../../bibliographie/ref}

\end{document}